\documentclass{article}
\usepackage{amsthm}
\usepackage{amssymb}
\usepackage{amsmath}
\usepackage{epsfig}
\usepackage{colonequals}
\usepackage{enumerate}
\usepackage{hyperref}

\newtheorem{theorem}{Theorem}[section]
\newtheorem{corollary}[theorem]{Corollary}
\newtheorem{lemma}[theorem]{Lemma}

\begin{document}
\title{Fine structure of $4$-critical triangle-free graphs II.  Planar triangle-free graphs with two precolored $4$-cycles}
\author{%
     Zden\v{e}k Dvo\v{r}\'ak\thanks{Computer Science Institute (CSI) of Charles University,
           Malostransk{\'e} n{\'a}m{\v e}st{\'\i} 25, 118 00 Prague, 
           Czech Republic. E-mail: \protect\href{mailto:rakdver@iuuk.mff.cuni.cz}{\protect\nolinkurl{rakdver@iuuk.mff.cuni.cz}}.
           Supported by project 14-19503S (Graph coloring and structure) of Czech Science Foundation.}
\and	   
Bernard Lidick\'y\thanks{
Iowa State University, Ames IA, USA. E-mail:
\protect\href{mailto:lidicky@iasate.edu}{\protect\nolinkurl{lidicky@iastate.edu}}.
Supported by NSF grant DMS-1266016.}
}
\date{\today}
\maketitle
\begin{abstract}
We study $3$-coloring properties of triangle-free planar graphs $G$ with two precolored $4$-cycles $C_1$ and $C_2$
that are far apart.  We prove that either every precoloring of $C_1\cup C_2$ extends to a $3$-coloring
of $G$, or $G$ contains one of two special substructures which uniquely determine which $3$-colorings of $C_1\cup C_2$ extend.
As a corollary, we prove that there
exists a constant $D>0$ such that if $H$ is a planar triangle-free graph and $S\subseteq V(H)$ consists
of vertices at pairwise distances at least $D$, then every precoloring of $S$ extends to a $3$-coloring of $H$.
This gives a positive answer to a conjecture of Dvo\v{r}\'ak, Kr\'al' and Thomas, and implies an exponential lower bound
on the number of $3$-colorings of triangle-free planar graphs of bounded maximum degree.
\end{abstract}

\section{Introduction}

The interest in the $3$-coloring properties of planar graphs was started by a celebrated theorem of Gr\"otzsch~\cite{grotzsch1959},
who proved that every planar triangle-free graph is $3$-colorable.  This result was later generalized and strengthened
in many different ways.  The one relevant to the topic of this paper concerns graphs embedded in surfaces.
While the direct analogue of Gr\"otzsch's theorem is false for any surface other than the sphere, $3$-colorability of
triangle-free graphs embedded in a fixed surface is nowadays quite well understood.

A first step in this direction was taken by Thomassen~\cite{thom-torus}, who proved that every graph of girth at least five
embedded in the projective plane or torus is $3$-colorable.  Thomas and Walls~\cite{tw-klein} extended this result
to graphs embedded in the Klein bottle.  More generally, Thomassen~\cite{thomassen-surf} proved that for any fixed surface $\Sigma$,
there are only finitely many $4$-critical graphs of girth at least $5$ that can be drawn in $\Sigma$ (a graph
$G$ is \emph{$4$-critical} if it is not $3$-colorable, but every proper subgraph of $G$ is $3$-colorable).  In other words,
there exists a constant $c_\Sigma$ such that if a graph of girth at least $5$ drawn in $\Sigma$ is not $3$-colorable,
then it contains a subgraph with at most $c_\Sigma$ vertices that is not $3$-colorable.  Thomassen's bound on $c_\Sigma$
is double exponential in the genus of $\Sigma$.  This was improved by Dvo\v{r}\'ak, Kr\'al and Thomas~\cite{trfree3},
who gave a bound on $c_\Sigma$ linear in the genus of $\Sigma$.  

In follow-up papers~\cite{trfree4,trfree6}, they
also described structure of $4$-critical triangle-free graphs on surfaces.
Essentially, they show that such a graph can be cut into a bounded number of pieces, each of them
satisfying one of the following: every $3$-coloring of the boundary of the piece $H$ extends to a $3$-coloring of $H$;
or a $3$-coloring of the boundary of $H$ extends to a $3$-coloring of $H$ if and only if it satisfies a specific condition
(``winding number constraint''); or the piece is homeomorphic to the cylinder and either all the faces of the piece have length $4$,
or both boundary cycles of the piece have length $4$ (the \emph{cylinder} is the sphere with two holes).
While the first two possibilities for the pieces of the structure give all the information about their $3$-colorings,
in the last subcase the information is much more limited.

The main aim of this series of papers is to fix this shortcoming.  
Let us remark that we cannot eliminate this subcase entirely---based on the construction of Thomas and Walls~\cite{tw-klein},
one can find $4$-critical triangle-free graphs $G_1$, $G_2$, \ldots drawn in a fixed surface $\Sigma$,
such that for any $i\ge 1$, $G_i$ contains two non-contractible $4$-faces bounding a subset $\Pi$ of $\Sigma$
homeomorphic to the cylinder, such that at least $i$ $5$-faces of $G$ are contained in $\Pi$.
Such a class of graphs cannot be described using only the first two subcases of the structure theorem.

Hence, the main problem studied in this paper concerns describing subgraphs of $4$-critical triangle-free graphs embedded in the cylinder
with rings of length $4$.
The exact description of such subgraphs under the additional assumption that all $(\le\!4)$-cycles are non-contractible was given by
Dvo\v{r}\'ak and Lidick\'y~\cite{dvolid}.
Even this special case is rather involved (in addition to one infinite class of such graphs mentioned in the previous paragraph,
there are more than 40 exceptional graphs) and extending it to the general triangle-free case would be difficult.
However, in the applications it is mostly sufficient to deal with the case that the boundary $4$-cycles of the cylinder are far apart,
and this is the case considered in this paper.  With this restriction, it turns out that there are only two infinite classes of critical graphs.

Before stating the result precisely, let us mention one application of this characterization.
\begin{lemma}\label{lemma-fourext}
There exists a constant $D\ge 0$ with the following property.
Let $G$ be a plane triangle-free graph, let $C$ be a $4$-cycle bounding a face of $G$
and let $v$ be a vertex of $G$.  Let $\psi$ be a $3$-coloring of $C+v$.
If the distance between $C$ and $v$ is at least $D$, then $\psi$ extends to a $3$-coloring of $G$.
\end{lemma}

This confirms Conjecture~1.5 of Dvo\v{r}\'ak et al.~\cite{trfree5}, who also proved that this 
implies the following more interesting result (stated as Conjecture~1.4 in~\cite{trfree5}).
\begin{corollary}\label{cor-farsv}
There exists a constant $D\ge 0$ with the following property.
Let $G$ be a planar triangle-free graph, let $S$ be a set of vertices of $G$ and let $\psi:S\to \{1,2,3\}$
be an arbitrary function.  If the distance between every two vertices of $S$ is
at least $D$, then $\psi$ extends to a $3$-coloring of $G$.
\end{corollary}

Thomassen~\cite{thom-many} conjectured that all triangle-free planar graphs have exponentially
many 3-colorings.  If $G$ is an $n$-vertex graph of maximum degree at most $\Delta$,
then there exists a set $S\subseteq V(G)$ of size at least $n/\Delta^D$ such that the distance between every two vertices of $S$ is
at least $D$.  Hence, Corollary~\ref{cor-farsv} implies that this conjecture
holds for triangle-free planar graphs of bounded maximum degree.

\begin{corollary}\label{cor-many}
Let $D$ be the constant of Corollary~\ref{cor-farsv}.
If $G$ is an $n$-vertex planar triangle-free graph of maximum degree at most $\Delta$,
then $G$ has at least $\left(3^{1/\Delta^D}\right)^n$ distinct $3$-colorings.
\end{corollary}

Let us now introduce a few definitions.  The \emph{disk} is the sphere with a hole.
In this paper, we generally consider graphs embedded in the sphere, the disk, or the cylinder.
Furthermore, we always assume that each face of the embedding contains at most one hole,
and if a face contains a hole, then the face is bounded by a cycle, which we call a \emph{ring}.
Note that the rings do not have to be disjoint.  Let $C$ be the union of the rings of such a graph $G$
($C$ is empty when $G$ is embedded in the sphere without holes).
We say that $G$ is \emph{critical} if $G\neq C$ and for every proper subgraph $G'$ of $G$ such that $C\subseteq G'$,
there exists a $3$-coloring of $C$ that extends to a $3$-coloring of $G'$, but not to a $3$-coloring of $G$; that is, removing any
edge or vertex not belonging to $C$ affects the set of precolorings of $C$ that extend to a $3$-coloring of the graph.

\begin{figure}
\begin{center}
\includegraphics[scale=0.8]{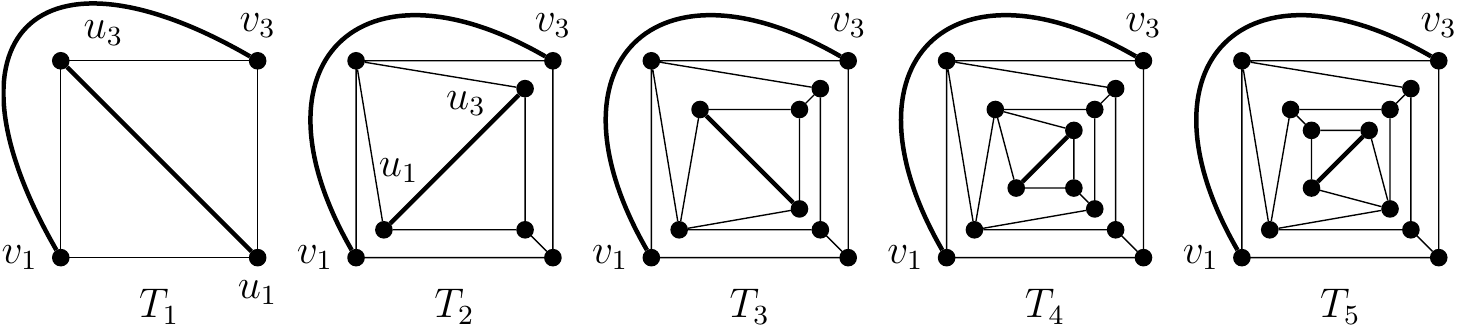}
\end{center}
\caption{Some Thomas-Walls graphs.}\label{fig-thomaswalls}
\end{figure}

We construct a sequence of graphs $T_1$, $T_2$, \ldots, which we call \emph{Thomas-Walls graphs} (Thomas and Walls~\cite{tw-klein} proved that they are exactly
the $4$-critical graphs that can be drawn in the Klein bottle without contractible cycles of length at most $4$).
Let $T_1$ be equal to $K_4$.  For $n\ge 1$, let $uv$ be any edge of $T_n$ that belongs to two triangles and let $T_{n+1}$ be obtained from $T_n-uv$
by adding vertices $x$, $y$ and $z$ and edges $ux$, $xy$, $xz$, $vy$, $vz$ and $yz$.  The first few graphs of this sequence are drawn in Figure~\ref{fig-thomaswalls}.
For $n\ge 2$, note that $T_n$ contains unique $4$-cycles $C_1=u_1u_2u_3u_4$ and $C_2=v_1v_2v_3v_4$ such that $u_1u_3,v_1v_3\in E(G)$.
Let $T'_n=T_n-\{u_1u_3,v_1v_3\}$.  We also define $T'_1$ to be a $4$-cycle.
We call the graphs $T'_1$, $T'_2$, \ldots \emph{reduced Thomas-Walls graphs}, and we say that $u_1u_3$ and $v_1v_3$ are their \emph{interface pairs}.
Note that $T'_n$ has an embedding in the cylinder with rings $C_1$ and $C_2$.  

A \emph{patch} is a graph drawn in the disk with ring $C$ of length $6$, such that $C$ has no chords and every face of the patch other than the one bounded by $C$ has length $4$.
Let $G$ be a graph embedded in the sphere, possibly with holes.
Let $G'$ be any graph which can be obtained from $G$ as follows.
Let $S$ be an independent set in $G$ such that every vertex of $S$ has degree $3$.  For each vertex $v\in S$
with neighbors $x$, $y$ and $z$, remove $v$, add new vertices $a$, $b$ and $c$ and a $6$-cycle
$C=xaybzc$, and draw any patch with ring $C$ in the disk bounded by $C$.
We say that any such graph $G'$ is obtained from $G$ by \emph{patching}.
This operation was introduced by Borodin et al.~\cite{4c4t} in the context of describing planar $4$-critical graphs with exactly $4$ triangles.

Consider a reduced Thomas-Walls graph $G=T'_n$ for some $n\ge 1$, with interface pairs $u_1u_3$ and $v_1v_3$.
A \emph{patched Thomas-Walls graph} is any graph obtained from such a graph $G$ by patching, and $u_1u_3$ and $v_1v_3$ are
its interface pairs (note that $u_1$, $u_3$, $v_1$, and $v_3$ have degree two in $G$, and thus they are not affected by patching).

Let $G$ be a graph embedded in the cylinder, with the rings $C_i=x_iy_iz_iw_i$
of length $4$, for $i=1,2$.  Let $y'_i$ be either a new vertex or $y_i$, and let $w'_i$ be either a new vertex
or $w_i$.  Let $G'$ be obtained from $G$ by adding $4$-cycles $x_iy'_iz_iw'_i$ forming the new rings.  We say that
$G'$ is obtained by \emph{framing on pairs $x_1z_1$ and $x_2z_2$}.

Let $G$ be a graph embedded in the cylinder with rings $C_1$ and $C_2$ of length $3$, such that every other face of $G$ has length $4$.
We say that such a graph $G$ is a \emph{$3,3$-quadrangulation}.
Let $G'$ be obtained from $G$ by subdividing at most one edge in each of $C_1$ and $C_2$.  We say that such a graph $G'$ is a \emph{near $3,3$-quadrangulation}.

We say that a graph $G$ embedded in the cylinder is \emph{tame} if $G$ contains no
contractible triangles, and all triangles of $G$ are pairwise vertex-disjoint.
The main result of this paper is the following.

\begin{theorem}\label{thm-mfar}
There exists a constant $D\ge 0$ such that the following holds.
Let $G$ be a tame graph embedded in the cylinder with rings $C_1$ and $C_2$ of length at most $4$.
If the distance between $C_1$ and $C_2$ is at least $D$ and $G$ is critical, then either
\begin{itemize}
\item $G$ is obtained from a patched Thomas-Walls graph by framing on its interface pairs, or
\item $G$ is a near $3,3$-quadrangulation.
\end{itemize}
\end{theorem}

This clearly implies the previously mentioned result on precolored $4$-cycle and vertex.

\begin{proof}[Proof of Lemma~\ref{lemma-fourext}]
Suppose for a contradiction that there exists a precoloring $\psi$ of $C+v$ that does not extend to a $3$-coloring of $G$.
Let $G'$ be obtained from $G$ by adding three new vertices $u_1$, $u_2$ and $u_3$ and edges $vu_1$, $u_1u_2$, $u_2u_3$ and $u_3v$,
with the $4$-cycle $C'=vu_1u_2u_3$ drawn so that it forms a face.  Drill a hole inside $C'$, turning it into a ring.
Let $\psi'$ be any $3$-coloring of $C\cup C'$ that extends $\psi$, and note that $\psi'$ does not extend to a $3$-coloring of $G'$.
Hence, $G'$ has a critical subgraph $G''$.
By Theorem~\ref{thm-many}, either $G''$ is obtained from a patched Thomas-Walls graph by framing on its interface pairs, or
$G''$ is a near $3,3$-quadrangulation.  However, this is a contradiction, since in either of the cases, each ring contains
at most two vertices of degree two, while $C'$ contains three vertices of degree two.
\end{proof}

Let us remark that which $3$-colorings of the rings of a framed patched Thomas-Walls graph or a near $3,3$-quadrangulation
extend to a $3$-coloring of the whole graph was determined in the first paper of this series (Dvo\v{r}\'ak and Lidick\'y~\cite{cylgen-part1}, Lemmas~2.5 and 2.9).

In Section~\ref{sec-cutting}, we prove a lemma based on an idea of~\cite{trfree5} that allows us to introduce new
non-contractible $(\le\!4)$-cycles.  This enables us to apply the result of the first paper of the series~\cite{cylgen-part1},
where we considered graphs in cylinder containing many non-contractible $(\le\!4)$-cycles, and thus prove Theorem~\ref{thm-mfar}
(Section~\ref{sec-mfar}).

\section{Cutting a cylinder}\label{sec-cutting}

Let us start by introducing several previous results.
The following is a special case of the main result of Dvo\v{r}\'ak et al.~\cite{trfree4}.

\begin{theorem}[\cite{trfree4}]\label{thm-critsize}
For every integer $k\ge 0$, there exists a constant $\beta$ with the following property.  Let $G$ be a graph embedded either in
the disk or the cylinder such that the sum of the lengths of the rings is at most $k$.
Suppose that $G$ contains no contractible triangles, and that every non-contractible $(\le\!4)$-cycle in $G$
is equal to one of the rings.  Let $F$ be the set of faces of $G$.  If $G$ is critical, then
\begin{align}
\sum_{f\in F} (|f|-4)\le \beta. \label{eq:beta}
\end{align}
\end{theorem}
Notice that \eqref{eq:beta} gives an upper bound $\beta$ on the number of $(\ge\!5)$-faces as well as upper bound $5\beta$  on number of vertices incident to $(\ge\!5)$-faces in $G$.

Furthermore, we need a result of Gimbel and Thomassen~\cite{gimbel}, which essentially states that in a plane triangle-free graph, every precoloring of a facial cycle $C$
of length at most $6$ extends to a $3$-coloring, unless $|C|=6$ and $G$ contains a subgraph such that every face other than the one bounded by $C$ has length $4$.

\begin{theorem}\label{thm-6cyc}
Let $G$ be a triangle-free graph drawn in the disk with the ring $C$ of length at most $6$.  If $G$ is critical, then
$|C|=6$ and all other faces of $G$ have length $4$.  Furthermore, a precoloring $\psi$ of $C=v_1v_2v_3v_4v_5v_6$
does not extend to a $3$-coloring of $G$ if and only if either
\begin{itemize}
\item $C$ has a chord $v_iv_{i+3}$ for some $i\in\{1,2,3\}$ and $\psi(v_i)=\psi(v_{i+3})$, or
\item $C$ is an induced cycle and $\psi(v_1)=\psi(v_4)$, $\psi(v_2)=\psi(v_5)$, and $\psi(v_3)=\psi(v_6)$.
\end{itemize}
\end{theorem}

Let us give an observation about critical graphs that is often useful.

\begin{lemma}\label{lemma-fr}
Let $G$ be a graph drawn in the sphere with holes so that every triangle is non-contractible.
If $G$ is critical, then
\begin{itemize}
\item every vertex $v\in V(G)$ that does not belong to the rings has degree at least three,
\item every contractible $(\le\!5)$-cycle in $G$ bounds a face, and
\item if $K$ is a closed walk in $G$ forming the boundary of an open disk $\Lambda$, then either $\Lambda$ is a face of $G$,
or all faces of $G$ contained in $\Lambda$ have length $4$.
\end{itemize}
\end{lemma}
\begin{proof}
Let $C$ be the union of the rings of $G$, and consider any vertex $v\in V(G)\setminus V(C)$.  Suppose for a contradiction
that $v$ has degree at most $2$.  Let $\psi$ be any $3$-coloring of $C$ that extends to a $3$-coloring $\varphi$ of $G-v$.
Then $\psi$ also extends to a $3$-coloring of $G$, by giving the vertices of $V(G)\setminus\{v\}$ the same color as in the
coloring $\varphi$ and by choosing a color of $v$ distinct from the colors of its neighbors.  This contradicts the assumption
that $G$ is critical.

The other two conclusions of the lemma follow by a similar argument using Theorem~\ref{thm-6cyc}.
\end{proof}

Let $G_1$ and $G_2$ be two graphs embedded in the cylinder, with the same rings.  We say that $G_1$ \emph{dominates} $G_2$ if every precoloring $\psi$ of the rings
that extends to a $3$-coloring of $G_1$ also extends to a $3$-coloring of $G_2$.
For two subgraphs $H_1$ and $H_2$ of a graph $G$, let $d(H_1, H_2)$ denote the length of the
shortest path in $G$ with one end in $H_1$ and the other end in $H_2$.  We now prove the key result of this section.

\begin{lemma}\label{lemma-cut}
There exists an integer $d_0\ge 3$ with the following property.
Let $G$ be a critical graph embedded in the cylinder with rings $C_1$ and $C_2$ of length at most $4$, such that
the distance between $C_1$ and $C_2$ is $d\ge d_0$, $G$ contains no contractible triangles, and every non-contractible $(\le\!4)$-cycle of $G$ is equal to one of the
rings.  There exists a tame graph $G'$ in the cylinder with the same rings, such that
\begin{itemize}
\item $G'$ dominates $G$,
\item the distance between $C_1$ and $C_2$ in $G'$ is at least $d-2$,
\item $G'$ contains a non-contractible $(\le\!4)$-cycle distinct from $C_1$ and $C_2$,
\item there exists a vertex $z\in V(G')$ such that the distance between $C_1\cup C_2$ and $z$ is at least three,
and $z$ is contained in all non-contractible $(\le\!4)$-cycles of $G'$ distinct from the rings, and
\item if $H'$ is a near $3,3$-quadrangulation with rings $C_1$ and $C_2$ and $H'\subseteq G'$, then $G$ is a near $3,3$-quadrangulation
with the same $5$-faces as $H'$.
\end{itemize}
\end{lemma}
\begin{proof}
Let $\beta$ be the constant of Theorem~\ref{thm-critsize}.  Let $d_0=(5\beta+1)(5(\beta+4)+2)+12$.
We prove the claim by induction on the number of vertices of $G$; hence, assume that the claim holds for all graphs with fewer than $|V(G)|$ vertices.

Let $f=v_1v_2v_3v_4$ be a $4$-face of $G$ at distance at least $3$ from $C_1\cup C_2$ such that no vertex of $f$ is incident with
a face of length greater than $4$.  Let $G_1$ be the graph obtained from $G$ by identifying $v_1$ with $v_3$ to a new vertex $z_1$, and let $G_2$ be the graph obtained from $G$ by
identifying $v_2$ with $v_4$ to a new vertex $z_2$.  If $G_1$ contained a contractible triangle $z_1xy$, then $G$ would contain a contractible $5$-cycle $v_1v_2v_3xy$.  Since
$G$ is critical, Lemma~\ref{lemma-fr} implies that $v_1v_2v_3xy$ bounds a face, contrary to the assumption that all faces incident
with $v_2$ have length $4$.
We conclude that every triangle in $G_1$ is non-contractible.  By the assumption on the distance between $f$ and $C_1\cup C_2$, the triangles of $G_1$ that contain $z_1$ are disjoint
with $C_1\cup C_2$.

Let us first discuss the case that the distance between $C_1$ and $C_2$ in $G_1$ is at least $d$.  
Observe that every $3$-coloring of $G_1$ extends to a $3$-coloring of $G$ (by giving $v_1$ and $v_3$ the
color of $z_1$), and thus $G_1$ dominates $G$.  Let $G'_1$ be a maximal critical subgraph of $G_1$, and note that $G'_1$ also dominates $G$.
If $G'_1$ is tame, then let $G''_1=G'_1$.  Otherwise, let $T_1$ and $T_2$ be triangles of
$G'_1$ containing $z_1$, chosen so that the part $\Sigma$ of the cylinder between $T_1$ and $T_2$ is maximal.
Let $G''_1$ be the graph obtained from $G'_1$ by removing all vertices and edges contained in the interior of $\Sigma$ and by identifying
the corresponding vertices of $T_1$ and $T_2$.

In the latter case, Theorem~\ref{thm-6cyc} implies that every $3$-coloring of $G''_1$ extends to a $3$-coloring of $G'_1$,
and thus $G''_1$ also dominates $G$.
Note that $G''_1$ is tame and all non-contractible $(\le\!4)$-cycles in $G''_1$ distinct from the rings contain $z_1$.

Let us consider the situation that a near $3,3$-quadrangulation $H'$ with rings $C_1$ and $C_2$ is a subgraph of $G''_1$.
Consider any face $h$ of $H'$.  Then $h$ either corresponds to a face in $G$, or $h$ is incident with $z_1$ and corresponds
to a contractible $(|h|+2)$-cycle $K$ in $G$ containing the path $v_1v_2v_3$.  Suppose the latter.
Since the distance between $z_1$ and $C_1\cup C_2$ in $H'$ is at least three, it follows that $h$ has length $4$.  Since $v_2$ is
not incident with a $6$-face in $G$, all the faces of $G$ contained in the closed disk bounded by $K$ have length $4$ by Lemma~\ref{lemma-fr}.
Since $v_1$ is not incident with a $6$-face in $G$, the same argument shows that if $G''_1\neq G'_1$, then all faces of $G$
contained in the disk bounded by the closed walk corresponding to $T_1\cup T_2$ have length $4$.  Note that one of these cases accounts
for all faces of $G$, and thus $G$ is a near $3,3$-quadrangulation with the same $5$-faces as $H'$.

If $G''_1$ contains a non-contractible $(\le\!4)$-cycle distinct from the rings, then the distance between the
rings of $G''_1$ is at least $d-2$ (due to the possible identification of $T_1$ with $T_2$ during the construction of $G''_1$)
and we can set $G'=G''_1$ in the conclusions of Lemma~\ref{lemma-cut}.  Otherwise, $G''_1=G'_1$, and thus $G''_1$ is critical and
the distance between the rings of $G''_1$ is at least $d$.  By the induction hypothesis, there exists a graph $G'$ dominating
$G'_1$ (and thus also $G$) that satisfies the conclusions of Lemma~\ref{lemma-cut}.

Therefore, we can assume that the distance between $C_1$ and $C_2$ in $G_1$ is less than $d$, and by symmetry, the distance between $C_1$ and $C_2$ in $G_2$ is also less than $d$.
Thus, without loss of generality, there exist paths $P_1$ and $P_2$ joining $v_1$ and $v_2$, respectively, to $C_1$, and paths $P_3$ and $P_4$ joining $v_3$ and $v_4$, respectively, to $C_2$,
such that $|P_1|+|P_3|\le d-1$ and $|P_2|+|P_4|\le d-1$.  Since the distance between $C_1$ and $C_2$ in $G$ is at least $d$, we have $|P_1|+|P_4|\ge d-1$ and $|P_2|+|P_3|\ge d-1$.
Summing the inequalities, we have $2d-2\le |P_1|+|P_2|+|P_3|+|P_4|\le 2d-2$, and thus all the inequalities hold with equality.  Consequently, $|P_1|=|P_2|$ and $|P_3|=|P_4|$.

Therefore, we can assume that for every $4$-face $f$ of $G$ at distance at least $3$ from $C_1\cup C_2$ that shares no vertex with a $(\ge\!5)$-face, there exists a labelling $v_1v_2v_3v_4$ of the vertices of $f$ such that $d(C_1,v_1)=d(C_1,v_2)$
and $d(C_1,v_3)=d(C_1,v_4)=d(C_1,v_1)+1$.  For an integer $a\ge 0$, let $A_a$ denote the set of vertices of $G$ at distance exactly $a$ from $C_1$.  Suppose that $5\le a\le d_0-5$
and that all vertices of $G$ in $A_{a-1}\cup A_a\cup A_{a+1}$ are only incident with $4$-faces.
Observe that if two vertices of $A_a$ are incident with the same face, then they are adjacent, and that $C_1$ and $C_2$ are in different components of $G-A_a$.  Therefore, there exists a non-contractible
cycle in $G$ with all vertices in $A_a$; let $Q_a$ denote the shortest such cycle (in particular, $Q_a$ is induced).

Let $Q_a=v_1\ldots v_k$.  For $1\le i\le k$, let $P_i$ denote a path of length $a$ from $v_i$ to $C_1$.  Observe that we can choose the paths so that for any $i,j\in \{1,\ldots, k\}$,
if $P_i$ and $P_j$ intersect, then $P_i\cap P_j$ is a path with an end in $C_1$.  Let $P_{k+1}=P_1$ and $v_{k+1}=v_1$.  For $1\le i\le k$, if $P_i$ and $P_{i+1}$ intersect, then let $K_i$ be the cycle contained
in $P_i\cup P_{i+1}\cup \{v_iv_{i+1}\}$.  Observe that since $C_1$ has at most $4$ edges, the cycle $K_i$ exists for at least $k-4$ values of $i$.  On the other hand, $K_i$ has odd length, and thus
an odd face of $G$ is contained in the open disk bounded by $K_i$.  Therefore, the cycle $K_i$ exists for at most $\beta$ values of $i$.  We conclude that $k\le\beta+4$.

By Theorem~\ref{thm-critsize}, at most $5\beta$ vertices of $G$ are incident with a face of length greater than $4$.
By the choice of $d_0$, there exists an integer $a\ge 5$ such that all vertices of $G$ at distance at least $a-1$ and at most $a+5(\beta+4)+1\le d_0-4$ from $C_1$ are incident only with $4$-faces.
Since $|Q_a|\le \beta+4$, there exists an integer $b$ such that $a\le b\le a+5(\beta+4)-5$ and $|Q_b|\le |Q_j|$ for all $b\le j\le b+5$.

Let $Q_b=v_1\ldots v_k$ and consider any $i\in\{1,\ldots,k\}$. Since no $4$-face has three vertices at distance $b$ from $C_1$, there exists an edge $v_iu_i\in E(G)$ such that $Q_b$ separates $u_i$ from $C_1$.
If there existed another such edge incident with $v_i$, then $G$ would contain a $4$-face $u_iv_iu'_iw$ with $d(u_i,C_1)=d(u'_i,C_1)=b+1$, which is a contradiction.
It follows that there exists exactly one such vertex $u_i$ for each $i=1,\ldots, k$.  Since all faces incident with $Q_b$
have length $4$, $u_1u_2\ldots u_k$ is a non-contractible closed walk in $G$.
By the choice of $Q_b$, $u_1u_2\ldots u_k$ is a cycle, and we can choose it as $Q_{b+1}$.

\begin{figure}
\begin{center}
\includegraphics[scale=1.0]{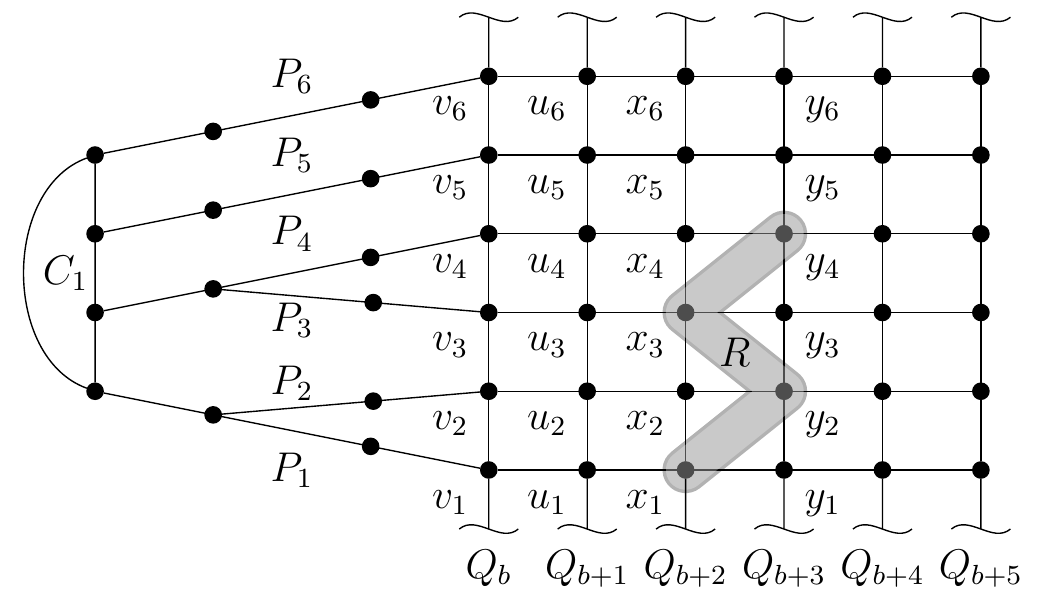}
\end{center}
\caption{Situation at the end of Lemma~\ref{lemma-cut} for $k=6$.}\label{fig-cut}
\end{figure}

The same argument shows that we can choose $Q_{b+2}$, \ldots, $Q_{b+5}$ so that the subgraph of $G$ drawn between $Q_b$ and $Q_{b+5}$ is a $6\times k$ cylindrical
grid, see Figure~\ref{fig-cut}.  Let $Q_{b+2}=x_1x_2\ldots x_k$ and $Q_{b+3}=y_1y_2\ldots y_k$.  If $k$ is even, then let $R=\{x_1,y_2,x_3,y_4,\ldots, x_{k-3}\}$.
If $k$ is odd, then let $R=\{x_1,y_2,x_3,\ldots, x_{k-2}\}$.  Let $G'$ be the graph obtained from $G$ by identifying the vertices of $R$ to a
single vertex $r$.  Observe that $G'$ dominates $G$ and that $r$ is contained in a unique non-contractible $(\le\!4)$-cycle in $G'$.
Furthermore, since we only contracted $4$-faces inside a quadrangulated cylinder, observe that if $G'$ contains a near $3,3$-quadrangulation
with rings $C_1$ and $C_2$ as a subgraph, then $G$ is a near $3,3$-quadrangulation as well.
Hence, $G'$ satisfies the conclusions of Lemma~\ref{lemma-cut}.
\end{proof}

We need an iterated version of this lemma.
Let $G$ be a tame graph embedded in the cylinder with rings of length at most $4$.
We say that $G$ is a \emph{chain} of graphs $G_1$, \ldots, $G_n$, if there exist non-contractible $(\le\!4)$-cycles $C_0$, \ldots, $C_n$
in $G$ such that 
\begin{itemize}
\item the cycles are pairwise vertex-disjoint except that for $(i,j) \in \{(0,1),(n,n-1)\}$,  $C_i$ can intersect $C_j$ if
$C_i$ is a 4-cycle and $C_j$ is a triangle,
\item for $0\le i<j<k\le n$ the cycle $C_j$ separates $C_i$ from $C_k$,
\item the cycles $C_0$ and $C_n$ are the rings of $G$,
\item every triangle of $G$ is equal to one of $C_0$, \ldots, $C_n$, and
\item for $1\le i\le n$, the subgraph of $G$ drawn between $C_{i-1}$ and $C_i$ is isomorphic to $G_i$.
\end{itemize}
We say that $C_0$, \ldots, $C_n$ are the \emph{cutting cycles} of the chain.

\begin{lemma}\label{lemma-manycuts}
For every integer $c>0$, there exists an integer $d_1\ge 0$ with the following property.
Let $G$ be a tame graph embedded in the cylinder with rings $C_1$ and $C_2$ of length at most $4$, such that
the distance between $C_1$ and $C_2$ is $d\ge d_1$.  If $G$ is not a chain of at least $c$ graphs, then
there exists a tame graph $G'$ in the cylinder with the same rings, such that
\begin{itemize}
\item $G'$ dominates $G$,
\item the distance between $C_1$ and $C_2$ in $G$ is at least $d-6c$,
\item $G'$ is a chain of $c$ graphs $G_1$, \ldots, $G_c$,
\item every non-contractible $(\le\!4)$-cycle in $G'$ intersects a ring of one of $G_1$, \ldots, $G_c$, and
\item if $G'$ contains a near $3,3$-quadrangulation with rings $C_1$ and $C_2$ as a subgraph, then
$G$ contains a near $3,3$-quadrangulation with rings $C_1$ and $C_2$ as a subgraph as well.
\end{itemize}
\end{lemma}
\begin{proof}
Let $d_0$ be the constant from Lemma~\ref{lemma-cut}, and let $d_1=c(d_0+7)+6c$.

Let $\mathcal{C}=K_0, K_1, \ldots, K_t$ be a sequence of cutting cycles of a chain in $G$ with $t<c$, chosen so that 
every non-contractible $(\le\!4)$-cycle in $G$ intersects one of the cutting cycles.
We prove that the conclusions of the lemma hold if the distance between $C_1$ and $C_2$
is at least $c(d_0+7)+6(c-t)$.  The proof is by induction on decreasing $t$; hence, we assume that the claim holds
for all chains of length greater than $t$.

Since the distance between $K_0$ and $K_t$ is greater than $c(d_0+7)$, there exists $i$ such that the distance between $K_i$ and
$K_{i+1}$ is at least $d_0+4$.  Let $K'_i$ and $K'_{i+1}$ be non-contractible $(\le\!4)$-cycles in $G$ intersecting $K_i$ and $K_{i+1}$,
respectively, such that the subgraph $F$ of $G$ drawn between $K'_i$ and $K'_{i+1}$ is minimal.  Let $F_0$ be a maximal critical
subgraph of $F$.  Note that the distance between the rings of $F_0$ is at least $d_0$, and that $F_0$ contains no non-contractible
$(\le\!4)$-cycles distinct from the rings by the choice of $\mathcal{C}$.  Let $F'$ be the graph obtained from $F_0$ by applying Lemma~\ref{lemma-cut}.
Let $Q$ be the graph obtained from $G$ by replacing the subgraph drawn between $K'_i$ and $K'_{i+1}$ by $F'$.

If $Q$ contains a near $3,3$-quadrangulation $H$ with rings $C_1$ and $C_2$ as a subgraph, then $H\cap F'$ is also
a near $3,3$-quadrangulation, and thus $F_0$ is a near $3,3$-quadrangulation with the same $5$-faces as $H\cap F'$
by the conclusions of Lemma~\ref{lemma-cut};
and thus the subgraph of $G$ obtained from $H$ by replacing $H\cap F'$ by $F_0$ is a near $3,3$-quadrangulation.

Since the distance between $K'_i$ and $K'_{i+1}$ in $F'$ is smaller than the distance in $F$ by at most two,
it follows that the distance between $C_1$ and $C_2$ in $Q$ is smaller than the distance in $G$ by at most $6$.
Let $\mathcal{C}'$ be the sequence obtained from $\mathcal{C}$ by adding one of the $(\le\!4)$-cycles of $F'$ distinct from its rings.
Observe that every non-contractible $(\le\!4)$-cycle in $Q$ intersects one of the cycles of $\mathcal{C}'$.

If $t=c-1$, then Lemma~\ref{lemma-manycuts} holds with $G'=Q$.  Otherwise, Lemma~\ref{lemma-manycuts} follows by the induction hypothesis for $Q$ and $\mathcal{C}'$.
\end{proof}

\section{Graphs in cylinder with rings far apart}\label{sec-mfar}

In \cite{cylgen-part1}, we proved the following.

\begin{theorem}[\cite{cylgen-part1}, Theorem 1.1]\label{thm-many}
There exists a constant $c\ge 0$ such that the following holds.
Let $G$ be a tame graph embedded in the cylinder with the rings $C_1$ and $C_2$ of length at most $4$.
If $G$ is a chain of at least $c$ graphs, then one of the following holds:
\begin{itemize}
\item every precoloring of $C_1\cup C_2$ extends to a $3$-coloring of $G$, or
\item $G$ contains a subgraph obtained from a patched Thomas-Walls graph by framing on its interface pairs, with rings $C_1$ and $C_2$, or
\item $G$ contains a near $3,3$-quadrangulation with rings $C_1$ and $C_2$ as a subgraph.
\end{itemize}
\end{theorem}

It is now straightforward to prove our main result by combining Theorem~\ref{thm-many} and Lemma~\ref{lemma-manycuts}.

\begin{proof}[Proof of Theorem~\ref{thm-mfar}]
Let $c$ be the constant of Theorem~\ref{thm-many}, and let $d_1$ be the constant of Lemma~\ref{lemma-manycuts} with
this $c$.  Let $D=\max(d_1, 15c)$.

Since $G$ is critical, there exists a precoloring of $C_1\cup C_2$ that does not extend to a $3$-coloring of $G$.
If $G$ is a chain of at least $c$ graphs, then the conclusions of Theorem~\ref{thm-mfar} follow from Theorem~\ref{thm-many},
since all contractible cycles of length at most $5$ in $G$ bound faces by Lemma~\ref{lemma-fr}.

Hence, assume that $G$ is not a chain of at least $c$ graphs, and let $G'$ be the graph obtained by Lemma~\ref{lemma-manycuts}.
Since $G'$ dominates $G$, there exists a $3$-coloring of $C_1\cup C_2$ that does not extend to a $3$-coloring of $G'$.
The graph $G'$ is a chain of graphs $G_1$, \ldots, $G_c$, such that every every non-contractible $(\le\!4)$-cycle in $G'$ intersects
one of their rings.  By Theorem~\ref{thm-many}, there exists a subgraph $H$ of $G'$ with rings $C_1$ and $C_2$ such that
either $H$ is obtained from a patched Thomas-Walls graph by framing on its interface pairs, or $H$ is a near $3,3$-quadrangulation.
In the latter case, the last condition in the conclusions of Lemma~\ref{lemma-manycuts} implies that $G$ is a near $3,3$-quadrangulation.

Hence, assume the former.  However, since the distance between $C_1$ and $C_2$ in $G'$ is at least $9c$, there exists $i\in\{1,\ldots,c\}$
such that the distance between the rings of $G_i$ is at least $5$.  But since $H\subseteq G$, it follows that $G_i$ contains
a non-contractible $4$-cycle disjoint from its rings, which is a contradiction.
\end{proof}

\bibliographystyle{acm}
\bibliography{cylgen}

\end{document}